\documentclass[11pt,a4paper]{article}
\usepackage{subfigure,amsmath,amssymb,amsfonts, amsmath}
\usepackage{a4wide}
\usepackage{verbatim}
\usepackage{epsfig,epic,eepic,graphicx}

\usepackage{amssymb}
\usepackage[latin1]{inputenc}
\usepackage{graphics}

\usepackage[francais]{babel}

\newtheorem{e-proposition}[theorem]{Proposition}

\newtheorem{e-definition}[theorem]{Definition\rm}

\newtheorem{theoreme}{Th\'eor\`eme}[section]

\newtheorem{corollaire}[theoreme]{Corollaire}

\newcommand{\al}{\alpha}
\newcommand{\bt}{\beta}

\newcommand{\be}{\begin{equation}}
\newcommand{\ee}{\end{equation}}

\newcommand{\bbR}{{\mathbb R}}

\def\og{\leavevmode\raise.3ex\hbox{$\scriptscriptstyle\langle\!\langle$~}}
\def\fg{\leavevmode\raise.3ex\hbox{~$\!\scriptscriptstyle\,\rangle\!\rangle$}}
\small

\begin{document}


\title{De Toda à KdV}



\author{D. Bambusi\footnote{Dipartimento di Matematica,                 
 Universit\`a degli Studi di Milano,          
 Via Saldini 50,                             
 20133 Milano,                              
 Italy, Dario.Bambusi@unimi.it}, T. Kappeler\footnote{Institut fuer Mathematik, Universitaet Zuerich,
Winterthurerstrasse 190, CH-8057 Zuerich, tk@math.unizh.ch}\hskip 0.2cm et T. Paul\footnote{CNRS et Département de Mathématiques et Applications UMR 8553 
\'Ecole Normale Supérieure, 45, rue d'Ulm, F 75730 Paris cedex 05, paul@dma.ens.fr}}


\date{}
\maketitle

\begin{abstract}
\selectlanguage{francais}
On considère la limite à grand nombre  de particules d'un 
système hamiltonien de type ``Toda p\'eriodique" pour  une famille de conditions initiales  proches de la solution d'équilibre. 
On montre  que, dans la formulation de paire de Lax, les deux bords des spectres des matrices de Jacobi des conditions initiales
sont déterminés, \`a une erreur pr\`es, par ceux de deux opérateurs de Hill, associ\'es \`a la famille de conditions initiales consid\'er\'ees. On en déduit que 
les spectres des matrices de Jacobi, lors de l'évolution limite donnée par KdV, restent  constants \` a une erreur pr\`es que nous estimons. 
Enfin on montre que les actions du syst\`eme Toda, convenablement renormalis\'ees, tendent vers celles des deux  \'equations de KdV.

\vskip 0.5\baselineskip
\end{abstract}

\section{Introduction et résultats}

Le système de Toda p\'eriodique à $N$ degrés de liberté est un système hamiltonien  où $N$ particules se meuvent sur le cercle réelle en 
interagissant avec un potentiel ``plus proche voisin" du type $e^{q_i-q_{i+1}}$ introduit dans \cite{T}. Il est bien connu, \cite{Fla}, que, par un changement (non symplectique) de coordonnées,
ce système peut se mettre sous la forme de ``Lax" $\dot L=[B,L]$, o\`u, apr\`es doublement de la dimension, 
\be\label{toda}\ 
L:=\left(\begin{array}{ccccc}
b_0&a_0&0&\dots&a_{2N-1}\\
a_{0}&b_1&a_{1}&\dots&0\\
0&a_{1}&b_2&\dots&0\\
0&\dots&\dots&\dots&a_{2N-2}\\
a_{2N-1}&\dots&\dots&a_{2N-2}&b_{2N-1}\end{array}\right)
 \ \mbox{,}\  
 \ee
 \be
B=\left(\begin{array}{ccccc}
0&a_0&0&\dots&-a_{2N-1}\\
-a_{0}&0&a_{1}&\dots&0\\
0&-a_{1}&0&\dots&0\\
0&\dots&\dots&\dots&a_{2N-2}\\
a_{2N-1}&\dots&\dots&-a_{2N-2}&0\end{array}\right),
\ee
et $a_i,b_i\in\bbR$ satisfont $a_{i+N}=a_i,\ b_{i+N}=b_i$. Nous allons consid\'erer des co\'efficients obtenus par discr\'etisation de fonctions lisses : 
$a_i=1+\epsilon^2\al(\frac i N)$ et 
$b_i=\epsilon^2\bt(\frac i{N})$,$\al,\ \bt\ \in C^\infty(\mathbb T)$ avec $|\epsilon|<<1$, une situation proche de la solution d'équilibre $a_i=1,\ b_i=0,\ i=0,\dots,2N-1$. 
De telles limites ont \'et\'e \'etudi\'ees par plusieurs auteurs (voir e.g. \cite{T}) dans des cas sp\'eciaux ; la nouveaut\'e pour le cas g\'en\'eral trait\'e ici est que, 
pour \'etudier cette limite, on a besoin de \textit{deux} op\'erateurs de Hill - voir aussi \cite{BP} et \cite{WS} o\`u se manisfeste ce ph\'enom\`ene dans l'\'etude de la dynamique
limite.

Rappelons que, si l'on considère l'opérateur de Schr\"odinger $H=-\frac{d^2}{dx^2}+u$ sur le cercle, alors l'\'equation de Korteveg-de Vries
\be
\label{kdv} \partial_tu-6u\partial_xu+ \partial^3_xu=0
\ee poss\`ede la forme de Lax  :
\be\label{laxkdv}
\dot H=[B,H]\ \ \ \mbox{où}\ \ \  B=-4\partial_x^3+6u\partial_x+3 u_x.\ee

D\'efinissons  $L_N=L^{\al,\bt}_N$, donn\'ee par (\ref{toda}) avec $a_i=1+\epsilon^2\al(\frac i N),\ b_i=\epsilon^2\bt(\frac i N)$ et $\epsilon=\frac 1 {2N}$.  
 Il est bien connu que les valeurs propres\  $\lambda^N_i,\ 0\leq i\leq 2N-1$ de $L_N$ sont toutes r\'eelles et vérifient
\begin{equation}
\label{ei.toda}
\lambda^N_0<\lambda^N_1\leq\lambda^N_2<...<\lambda^N_{2N-3}\leq\lambda_{2N-2}^N<\lambda_{2N-1}^N.  
\end{equation} 
En particulier, pour $\alpha=\beta=0$ (solution d'équilibre), le spectre de $L^{0,0}_N$ est donn\'e par :
\be\label{equi}
\lambda_0=-2,\quad \lambda_{2l-1}=\lambda_{2l}=-2\cos{\frac{l\pi}{N}}\ ,\quad l=1,\dots, N-1,\quad \lambda_{2N-1}=2.
\ee
Le c\oe ur de nos résultats est le Théorème suivant.
\begin{theoreme}\label{spectre}
Soit $L_N$ la matrice $L^{\al,\bt}_N$ avec $\al,\bt\in C^\infty(\mathbb T)$, $\int_{\mathbb T}\al(x)dx=\int_{\mathbb T}\bt(x)dx=0$ 
et soit $M_N=[N^{\frac 1 4}]$. Alors le spectre de $L_N$ 
est constitu\'e de $2N$ nombres r\'eels $(\lambda_i^N)_{i=0,\dots,2N-1}$ qui satisfont, pour tout $\delta>0$, uniformément 
par rapport à 
 $j=0,\dots, 2N-1$ et $\al,\ \bt$ dans tout ensemble born\'e de $C^\infty(\mathbb T)$,
\begin{eqnarray}
\label{thierry.1}
\lambda_j^N&=&-2+\frac 1 {4N^2}\lambda_j^-+O(N^{-3+\delta})\ ,\quad j=0,...,2M_N
\\
\label{thierry.2}
\lambda_{2l-1}^N, \lambda_{2l}^N&=&-2\cos{\frac{l\pi}{N}}+O(N^{-3+\delta})\
,\quad l=M_N+1,...,N-1-M_N 
\\
\label{thierry.3}
\lambda_{2N-1-j}^N&=&2-\frac 1 {4N^2}\lambda_j^++O(N^{-3+\delta})\ ,\quad j=0,...,2M_N
\end{eqnarray}
où $(\lambda^\pm_i)_{i=0,\dots,2M_N}$ sont les $2M_N+1$ premières valeurs propres des opérateurs de Hill :
\be\label{schrokdv}
H_\pm=-\frac{d^2}{dx^2}-2\al(x)\mp \bt(x).
\ee
\end{theoreme}
\begin{corollaire}\label{cor}
Soient $\al_t=-(u^-_t+u^+_t)/4$ et $\bt_t=(u^-_t-u^+_t)/2$ obtenues en faisant évoluer les conditions initiales $u^+=-2\al-\bt$ et $u^-=-2\al+\bt$ par l'équation de KdV. 
Soit $L^t_N=L_N^{\al_t,\bt_t}$ donnée par (\ref{toda}) et soit 
$\{\lambda^N_j(L^t_N)\}_{j=0,\dots,2N-1}$ son spectre (qui est  conservé par la dynamique de Toda). Alors, pour tout $\delta>0$ et uniformément en temps:
\be
\lambda^N_j(L^t_N)-\lambda^N_j(L^0_N)=O(N^{-3+\delta}).
\ee
\end{corollaire}

Le
 deuxième résultat concerne l'asymptotisme \`a grand $N$ des actions de Toda.
Rappelons, \cite{HK1}, que, pour $1\leq n\leq N-1$, la $n$i\`eme action d'une matrice de Jacobi $L_N$  v\'erifie  la formule 
\be\label{fond}I^N_n=\frac 1 \pi\int_{\lambda^N_{2n-1}}^{\lambda^N_{2n}}arcosh\left((-1)^{N-n}\frac{\Delta^N(\lambda)}2\right)d\lambda,
\ee 
o\`u 
$\Delta^N(\lambda)$ est le discriminant de $L_N$. 
On a  une formule similaire pour les actions $I^\pm_n$ des deux op\'erateurs de Hill $H_\pm$ d\'efinis par (\ref{schrokdv}) ; voir (\ref{hill}) plus bas.
\vskip 0.2cm
\begin{theoreme}\label{act} Pour tout $n\geq 1$,
\[8N^2I^N_n\to I^-_n\ \ \  \mbox{ et }\ \ \ \ 8N^2I^N_{N-n}\to I^+_n\ \ \ \  \mbox{ quand }\  N\to\infty.\]
\end{theoreme}
D'une fa\c con analogue, on obtient le comportement asymptotique \`a grand $N$ des fr\'equences du syst\`eme de  Toda \cite{BKP}.

%


\section{Quelques éléments de preuve}
\subsection{Lax et T\"oplitz}
L'id\'ee principale de la preuve du Th\'eor\`eme \ref{spectre} consiste \`a consid\'erer  la matrice de Jacobi $L_N$ comme la matrice dans une base canonique 
d'un op\'erateur de
T\"oplitz pour la quantification du tore $\mathbb T^2=\bbR^2/\mathbb Z^2$ avec une constante de Planck $\hbar=\frac1 {4\pi N}=\frac\epsilon{2\pi}$. Cette identification a \'et\'e propos\'ee par Bloch et al \cite{BGPU} pour \'etudier la limite \`a grand $N$ du 
syst\`eme de Toda p\'eriodique. Nous allons, dans cet article, nous placer dans
ce formalisme pour \'etudier le spectre de $L_N^{\al,\bt}$. 

Pour cela on considère l'espace de Hilbert de dimension $2N$, $\mathcal H_{2N}$, généré par les fonctions \textit{Theta} définies par :

\be\label{teta}
\Theta_j(z=x+iy):=(4N)^{1/4}e^{-\pi j^2/2N}\sum_{n\in\mathbb Z} e^{-\pi(2Nn^2+2jn)}e^{2\pi iz(j+2Nn)},\ee
$\ j=0,\dots,2N-1$.

Le produit scalaire $\langle.,.\rangle$ est celui de $L^2([0,1]\times[0,1],e^{-4\pi Ny^2}dxdy)$, de mani\`ere \`a ce que\ $\{\Theta_j\}_{j=0,\dots,2N-1}$ soit une base
orthonorm\'ee de $\mathcal H_{2N}$ . 
La matrice $L_N^{\al,\bt}$ est la matrice, dans la base $\{\Theta_{2N-1},\dots,\Theta_0\}$,
d'un opérateur de T\"oplitz $T^{\al,\bt}_N$ de symbole principal  (voir \cite{BGPU})$$\epsilon^2\bt(x)+2(1+\epsilon^2\al(x))\cos{2\pi y}.$$
Il est donc \'el\'ementaire de constater   que les vecteurs $\psi^k\in\mathcal H_{2N},\  k=0,\dots,2N-1,$ définis par:
\be\label{dsa}
\psi^k(z)
=\frac 1 {(2N)^{1/2}}\sum_{j=0}^{2N-1} e^{\pi i\frac{kj}N} \Theta_j(z)
\ee sont  vecteurs propres de l'op\'erateur $T^{0,0}_N$ dans $\mathcal H_{2N}$, 
 de valeurs
propres correspondantes  données par (\ref{equi}). On remarque  que 
\[
\psi^k(z)=(4N)^{-1/4}
\int_0^1\rho(z,k/2N+is)e^{-2\pi Ns^2}ds,\] o\`u  $\rho(z,\overline{z'}):=\sum_{j=0}^{2N-1}\Theta_j(z)\overline{\Theta_j(z')}$.
\subsection{Quasimodes}
Dans l'idée de \cite{PU} nous allons construire des quasimodes de $T_N^{\al,\bt}$  comme superpositions pondérées d'états cohérents $\rho(z,k/2N+is)$ et nous verrons 
que l'équation de  KdV va apparaître comme \textit{équation de transport}. 

Définissons, pour $\mu\in C^\infty(\mathbb T),\ \mu(x)=\sum_{k\in\mathbb Z}\mu_ke^{2\pi ikx}$, 

\be\label{wkb}
\psi^k_\mu(z)=(4N)^{-1/4}\int_0^1\rho(z,k/2N+is)\mu(s)e^{-2\pi Ns^2}ds.
\ee
Le résultat suivant, cons\'equence du calcul symbolique ``à la T\"oplitz" mais que l'on peut obtenir ici par un calcul direct, est le c\oe ur de la preuve :
\begin{theoreme}
Pour tout $\mu,\mu'\in C^\infty(\mathbb T)$ arbitraire,
\[
 <\psi^k_\mu,\psi^{k'}_{\mu'}>=\sum_{l\in\mathbb Z} \overline{\mu_l}\mu'_{l-k+k'}e^{-\pi l^2/2N}e^{-\pi (l-k+k')^2/2N}.\ 
 \]
  \[\mbox{En particulier \ \ \ \ }\ \ \ ||\psi^k_\mu||^2=\sum_{l\in\mathbb Z}|\mu_l|^2e^{-\pi l^2/N}\sim ||\mu||^2_{L^2(\mathbb T)}\ \mbox{quand}\ N\to+\infty.\]
 
De plus, si $T^{\al,\bt}_N$ est l'opérateur dont la matrice sur la base $\{\Theta_{2N-1},\dots,\Theta_0\}$ est $L_N^{\al,\bt}$ et $\epsilon=1/2N$, 
 \be
 T^{\al,\bt}_N\psi^k_\mu=\psi^k_{\mu^k}+O(\epsilon^{3})
 \ee
avec
\[
 \mu^k(x)=\left(-2\cos{(2\pi  k \epsilon- i \epsilon \frac d {dx})}+\epsilon^2\left(-2\al(x)\cos{(2\pi  k \epsilon- i \epsilon \frac d {dx})}+\bt(x)\right)\right)\mu(x).
 \]
 \end{theoreme}
On voit donc que $\psi^k_\mu$ sera un quasimode de $T^{\al,\bt}_N$ à l'ordre $\epsilon^3$ si $\mu$ est lui-même un quasimode de l'opérateur 
\[-2\cos{(2\pi  k \epsilon- i \epsilon \frac d {dx})}+\epsilon^2\left(-2\al(x)\cos{(2\pi  k \epsilon- i \epsilon \frac d {dx})}+\bt(x)\right).\]
 Plusieurs comportements sont à envisager:

- milieu du spectre :
 pour $l=M_N+1,\dots,N-1-M_N$ on a, pour les valeurs propres distinctes de $L_N^{0,0}$,  toutes de multiplicit\'e deux,
\[
|\cos{\pi (l+1)\epsilon}-\cos{\pi l\epsilon}|\sim 4  \pi   (2M_N+1)\epsilon^2\geq \epsilon^{7/4}>\epsilon^2 \  \mbox{ quand } \epsilon\to 0,
\]
et donc la théorie des perturbations (avec dégénérescence) s'applique. On montre facilement que la condition de moyenne 
nulle de $\al,\bt$ donne une  correction au premier ordre qui est en fait d'ordre $\epsilon^\infty$.

-  les bords du spectre, pour $0\leq k\leq 2M_N$ et $2N-1-2M_N\leq k\leq 2N-1$.

L'équation aux valeurs propres\[
\left(-2\cos{(2\pi  k \epsilon- i \epsilon \frac d {dx})}+\epsilon^2\left(-2\al(x)\cos{(2\pi  k \epsilon- i \epsilon \frac d {dx})}+\bt(x)\right)\right)\mu(x)
 =\lambda\mu(x)
 \]
 devient, avec l'ansatz $\mu(x)=e^{-2i\pi kx}\nu_k(x)$ (resp. $e^{-2i\pi (k-N)x}\nu_k(x))$ pour $k$ proche de $0$ (resp. $2N-1$),
 \[\epsilon^2\left(-\frac {d^2} {dx^2}- 2\al(x)\pm\bt(x)\right)\nu_k=(2\pm\lambda_k^N)\nu_k+O(\epsilon^3)\] suivant  le bord consid\'er\'e. On reconnaît là 
 les deux opérateurs de Hill 
 $H_\pm$. Un  argument de comptage montre qu'on obtient bien ainsi tout le spectre de $L_N^{\al,\bt}$, les ph\'enom\`enes de 
 d\'eg\'en\'erescence \'etant responsables du terme suppl\'ementaire en $N^\delta$  et la valeur de $M_N=[N^{1/4}]$ assurant 
 l'uniformit\'e de la transition entre les trois parties du spectre. Le Corollaire \ref{cor} se d\'emontre gr\^ace \`a l'uniformit\'e des estimations semiclassiques par
 rapport aux normes $sup$ des d\'eriv\'ees de symboles, et 
  au fait que les normes $sup$ restent born\'es uniform\'ement en temps lors de l'\'evolution par KdV \cite{BKP}.

\subsection{Intégrabilité}
La preuve du Th\'eor\`eme \ref{act} repose sur le fait que l'on puisse \'ecrire les discriminants $\Delta^N(\lambda)$ de la matrice de Jacobi $L_N$ et ceux, 
$\Delta^\pm(\lambda)$, de $H_\pm$ comme :
\[
\Delta^N(\lambda)^2-4=\prod_{j=0}^{2N-1}(\lambda_j^N-\lambda)\ \ \mbox{ et }\ \Delta^\pm(\lambda)^2-4=4\prod_{j\geq 0}\frac{\lambda^\pm_j-\lambda}{\pi_j^2},\] 
\[\pi_j:=\pi[(j+1)/2]\ (j\geq 1),\ \pi_0=1.\]

Soit $[\Lambda_1,\Lambda_2]$ un intervale compact de $\bbR$ avec $\Lambda_1\leq 0<\Lambda_2$.
Le Th\'eor\`eme \ref{act} est alors une cons\'equence du
\begin{theoreme}
\label{discriminant}
 Uniform\'ement pour $\Lambda_1\leq\lambda\leq\Lambda_2$, 
\begin{equation}
\label{dis.e}
\lim_{N\to\infty}(-1)^N\Delta^N(-2+\epsilon^2
  \lambda)=
  \Delta^-(\lambda)\ \mbox{ et }\ 
\lim_{N\to\infty}\Delta^N(2-\epsilon^2
  \lambda)=
  \Delta^+(\lambda),
\end{equation}
\end{theoreme}
et de la formule (\ref{fond}) pour les actions de Toda, ainsi que la formule suivante pour les variables actions $I^\pm_n$  de $H_\pm$:
\be\label{hill}
I^\pm_n=\frac 2 \pi\int_{\lambda^\pm_{2n-1}}^{\lambda^\pm_{2n}}arcosh\left((-1)^n\frac{\Delta^\pm(\lambda)}2\right)d\lambda.
\ee 





\end{document}